\documentstyle[12pt]{article}
\def\baselinestretch{1.0}
\title{\bf  Noise and chaotic disturbance \\ 
on self-similar set\footnotemark
}
\author{  Zu-Guo Yu \\
 {\small (Institute of Theoretical Physics, Academia Sinica},\\
  {\small P.O. Box 2735, Beijing 100080, P. R. China. yuzg@itp.ac.cn)}\\
    Fu-Yao Ren\\
  {\small (Institute of Mathematics, Fudan University, Shanghai 200433 P.R.C.)}
 }
\date{}
\begin{document}
\newcommand{\be}{\begin{equation}}
\newcommand{\ee}{\end{equation}}
\newtheorem{Theorem}{\quad Theorem}
\newtheorem{Proposition}{\quad Proposition}[section]
\newtheorem{Definition}{\quad Definition}[section]
\newtheorem{Lemma}{\quad Lemma}[section]
\newtheorem{Corollary}{\quad Corollary}[section]
\newtheorem{Example}{\quad Example}[section]
\maketitle
 \renewcommand{\thefootnote}{\fnsymbol{footnote}}
 \footnotetext{* Project partially supported by the Tianyuan Foundation of China.}
\begin{abstract}
The effect of noise on self-similar set is studied. The iteratie procedure
used to generate the self-similar set is moidified by adding a stochastic
variable to the diameter of generating sets at each iteration.  The noise may causes
the generating set to collapse
to a point. Distribution functions are found describing the probability that
any generating set collapse. The effect of chaotic
disturbance on the iteration of self-similar set is studied. It is shown that
the iterative procedure which describes the self-similar set is truncated under
the influence of disturbance generated by the tent map. Conditions which lead to
truncation of any chaotic map are also obtained. \end{abstract}

\section{Introduction}

\ \ \ \   Many dynamical systems have an attractor or a repellor that has in some way
a self-similar structure. In reality, fractal growth is always under the influence
of the environment. In order to simulate realistic fractal growth, it is important
to in corporate disturbances or noises into the iteration procedure.

  The effect of noise on chaotic dynamical systems is of great interest and has
been studied by many authors. The early work on this problem was carried out by
Crutchfield et al$^{\cite{CFH}}$, who studied the effect of noise on period doubling in
a discrete system. Additional work was carried out by Svensmark and 
Samuelson$^{\cite{SS}}$
on the Josephson junction. Wiesenfeld and McNamara$^{\cite{WM}}$ have studied the amplification
of a small resonant periodic perturbation in the presence of noise near the period doubling
threshold. Arecchi et al$^{\cite{ABP}}$ have studied the effect of noise on the forced
Duffing oscillator in the region of parameter space, where different chaotic attractors
coexist. Kautz$^{\cite{Kau}}$ has investigated the problem of thermally induced escape
from the basin of attraction in a dc-biased Josephson junction. Last, 
Kapi-taniak$^{\cite{Kapi}}$
has studied the behavior of the probability density function of a driven nonlinear
system, his result implies that the noise may introduce a degree of order in a
chaotic system, and the exponent is a random number and has a corresponding probability
density function.

  Some authors have studied the effects of noise on discrete dynamical systems.
For example, Crutchfield and Packard$^{\cite{CP}}$ have studied the symbolic dynamics
of chaotic maps when they are perturbed by a noise term. Carlson and 
Schieve$^{\cite{CS}}$
 considered
noise to the standard shift map. Garc\'{i}a-Pelayo and schieve$^{\cite{GPS}}$ introduced noise
to the affine contractive iterated function system. Cole and Schieve$^{\cite{CoS}}$ studied
effect noise on triadic Cantor set. More resently, Chia-Chu Chen$^{\cite{Chen}}$ studied
the effect of chaotic disturbance on triadic Cantor set. Chia-Chu Chen$^{\cite{Chen}}$ pointed
out that it will be more interesting if we can show that truncation of a fractal
set under the influence of chaotic noise is a general phenomenon. Based the thought
of ref.\cite{CoS,Chen}, in this paper we study the effect of noise and chaotic disturbance
on a class of more general fractal set, i.e. self-similar set.

\section{ Fractal constructure of self-similar  sets and noise}

  Denote the $d$-dimensional Euclidean space by ${\bf R}^d$, and
   fixed $K<\infty$. Let 
 and  $$S_j(x)=\xi_j R_jx+b_j,\ 0<\xi_j<1,\ R_j\ \mbox{ orthogonal,}\  b_j
 \in {\bf R}^d,\quad j\in J=\{1,2,\cdots K\} $$
be contractive similarities map on a compact set $E_0$ which satisfies $\overline{Int(E_0)}=E_0$,
  and assume that  $Int(S_j(E_0))\cap Int(S_i(E_0))=\emptyset$
for $j\neq i$.
  For natural number $n$, let
\begin{eqnarray*}
 E_{j_1j_2\cdots j_n}&=&S_{j_1}\circ \cdots \circ  S_{j_n}(E_0)\\
  E(n)=& &\cup_{j_i\in I}E_{j_1\cdots j_n}. 
  \end{eqnarray*}
It is obvious that
$$E_{j_1\cdots j_n}\subset E_{j_1\cdots j_{n-1}}, \quad E(n)\subset E(n-1).$$
Then $$E=\cap_{n=1}^{\infty}E(n)$$ is called {\it self-similar set}, and $E_{j_1\cdots j_n}$ is
called {\it generating set} of $n$th stage.
  From ref.\cite{Hutch},  the fractal dimension $s$ of $E$ is the solution of equation $\sum_{j=1}^{K}\xi_j^s=1.$

 {\it Examples}: Cantor's set , Cantor's $k$-bars,  Von Koch snowflake and Sierpinski gasket
 are
self-similar set.

  We might as well assume $|E_0|=1$, where $|E_0|$ denote the diameter of $E_0$. It is easy
  to see 
  \be L_0=|E_0|=1,\quad L_{j_1\cdots j_n}=|E_{j_1\cdots j_n}|=\xi_{j_1}
  \xi_{j_2}\cdots \xi_{j_n}. \label{nfs1} \ee
 The noise is introduced by the above rules with the addition of $K$ independent stochastic variables,
 $\delta_1, \cdots, \delta_K$. Under this correction, we have
 \be L_0=|E_0|=1,\quad L_{j_1\cdots j_n}=\xi_{j_n}(L_{j_1\cdots j_{n-1}} +
 \delta_{j_n}).\label{nfs2} \ee

  The generating set $E_{j_1\cdots j_n}$ is said to {\it 'collapse'} when $L_{j_1\cdots j_n}$
  becomes less than or equal to zero.
   We terminate the iteration of a
   generating set only when it  collapse.
But we may not terminate the iteration other generating set which is not subset of this
 generating set  as the same time.

\section{Distribution function}

 We want to obtain 
 \begin{Theorem}
 The  distribution function that describes the probability
  that the generating set will collapse on $n$th iteration for two different
   cases can be given out definitely.

For case 1, our 'noises' 
   $\delta_i$, $( i=1,2,\cdots,K)$ are stochastic variables
that can take the values $-\triangle_i$, $0$, $\triangle_i$ with a probability of $1/3$,
where $\triangle\in (0,1)$.
For case 2, our 'noises' $\delta_i$, $( i=1,2,\cdots,K)$ are stochastic variables
with an arbitrary normalized probability functions.
\end{Theorem}

  For any generating set $E_{j_1\cdots j_n}$, from (\ref{nfs2}) we have
 \be L_{j_1\cdots j_n}=\xi_{j_1}\cdots \xi_{j_n}+N_{j_1\cdots j_n},\label{nfs3}\ee
where
\be N_{j_1\cdots j_n}=\xi_{j_n}\delta_{j_n}+\xi_{j_n}\xi_{j_{n-1}}\delta_{j_{n-1}}+\cdots+
\xi_{j_n}\cdots \xi_{j_1}\delta_{j_1} \label{nfs4} \ee
denote the noise term. From (\ref{nfs2}) we also have
\be N_{j_1\cdots j_n}=\xi_{j_n}(N_{j_1\cdots j_{n-1}}+\delta_{j_n}).
 \label{nfs5} \ee
 From (\ref{nfs4}), the generating set collapse if and only if $N_{j_1\cdots j_n}\leq -\xi_{j_n}\cdots\xi_{j_1}$.

  We denote $C_{j_1\cdots j_n}$ the probability distribution describing the chance
  that the generating set $E_{j_1\cdots j_n}$ collapse, $NT_{j_1\cdots j_n}$
  the probability that the previous $n-1$ generating  sets $\{E_{j_1},\cdots,E_{j_1\cdots j_{n-1}}\}$
  do not collapse,$LE_{j_1\cdots j_n}$ the probability that the noise term is
  $\leq -\xi_{j_n}\cdots \xi_{j_1}$, $GE_{j_1\cdots j_n}$ the probability that the noise term is
  $\geq \xi_{j_n}\cdots \xi_{j_1}$. Then
 \be C_{j_1\cdots j_n}=NT_{j_1\cdots j_n}\cdot LE_{j_1\cdots j_n}. 
 \label{nfs6} \ee
It is easy to see $GE_{j_1\cdots j_n}=LE_{j_1\cdots j_n}$.

  \subsection{ Distribution function for case 1.}

  First, we determine $LE_{j_1\cdots j_n}$. 

   We denote $\xi=\max_{1\leq i \leq K}\xi_i$, $\triangle=\max_{1\leq i \leq K}\triangle_i$,
   $\triangle'=\min_{1\leq i \leq K} \triangle_i$. Since $\xi<1$, then
   $\sum_{n=1}^{\infty}\xi^n=\frac{\xi}{1-\xi}<\infty$. For $\xi$ and $\triangle,\ \triangle'$, we assume
  \be \xi\leq \frac{\triangle'}{2\triangle+\triangle'}. \label{nfs7} \ee
   From (\ref{nfs4}), we have
 \be -\frac{\xi\triangle}{1-\xi}\leq N_{j_1\cdots j_n}\leq 
 \frac{\xi\triangle}{1-\xi}. \label{nfs8} \ee

  For any $n$, we must have the following three cases.

 If $\xi_{j_1}\cdots\xi_{j_n}>\frac{\xi\triangle}{1-\xi}$, then from (\ref{nfs8}), the generating set $N_{j_1\cdots j_n}$
 can not $\leq -\xi_{j_1}\cdots\xi_{j_n}$, hence $LE_{j_1\cdots j_n}=0$.

 If $\xi_{j_1}\cdots\xi_{j_{n-1}}>\frac{\xi\triangle}{1-\xi}$ and $\xi_{j_1}\cdots\xi_{j_n}\leq \frac{\xi\triangle}{1-\xi}$,
since the possible values of $N_{j_1\cdots j_n} $ are evenly spaced, from ref.[10], the number
of these values in a given range is proportional to the length of the range,
the points of the numerator are confined to the region $[-\frac{\xi\triangle}{1-\xi},-\xi_{j_1}\cdots\xi_{j_n}]$,
while the points of the denominator are in $[-\frac{\xi\triangle}{1-\xi},
\frac{\xi\triangle}{1-\xi}]$, hence
\begin{eqnarray*}
LE   &=\frac{ \mbox{ number\ of\ possible\ values\ of}\ N_{j_1\cdots j_n}\ 
\mbox{ less\ than}\ -\xi_{j_1}\cdots\xi_{j_n} }
{ \mbox{ total\ number\ of\ possible\ values\ of}\ N_{j_1\cdots j_n} }\\
 &\simeq \frac{\frac{\xi\triangle}{1-\xi}-\xi_{j_1}\cdots
  \xi_{j_n}}{2\frac{\xi\triangle}{1-\xi}}=\frac{\xi\triangle-(1-\xi)\xi_{j_1}
  \cdots\xi_{j_n}}{2\xi\triangle}. 
  \end{eqnarray*}
 Our approximation becomes very good for $\triangle<<1$.

If $\xi_{j_1}\cdots\xi_{j_{n-1}}\leq \frac{\xi\triangle}{1-\xi}$, we have $-\xi_{j_1}\cdots \xi_{j_n}\geq -\xi_{j_n}\frac{\xi\triangle}{1-\xi}$.
 when $\delta_{j_n}=-\triangle_{j_n}$,
from (\ref{nfs5}) and (\ref{nfs8}), we have $\xi_{j_n}(-\frac{\xi\triangle}{1-\xi}-\triangle_{j_n})\leq N_{j_1\cdots j_n}\leq \xi_{j_n}(\frac{\xi\triangle}{1-\xi}-\triangle_{j_n})$.
From (\ref{nfs7}), we have $\xi_{j_n}(\frac{\xi\triangle}{1-\xi}-\triangle_{j_n})\leq -\xi_{j_n}\frac{\xi\triangle}{1-\xi} $.
Hence $N_{j_1\cdots j_n}\leq -\xi_{j_1}\cdots\xi_{j_n}$.
 This means that if the last step was negative then the generating set must collapse.
 But when $\delta_{j_n}\geq 0$,
since $N_{j_1\cdots j_{n-1}} >-\xi_{j_1}\cdots\xi_{j_{n-1}}$ and
$$N_{j_1\cdots j_n} =\xi_{j_n}N_{j_1\cdots j_{n-1}} +\xi_{j_n}\delta_{j_n},$$
we have $N_{j_1\cdots j_n}>-\xi_{j_1}\cdots\xi_{j_n}$. Hence the generating set
can collapse only if $\delta_{j_n}$ was negative. Hence
$$LE_{j_1\cdots j_n}=(\mbox{ probability\ that}\ \delta_{j_n}=-\triangle_{j_n})=1/3.$$

 Second, we determine $NT_{j_1\cdots j_n}$.

  We find that
 \be NT_{j_1\cdots j_n}=1-(\mbox{ probability\ that\ one\ of\ generating\ 
 sets}\ E_{j_1},\ \cdots,\ E_{j_{n-1}}\ \mbox{ collapse}),\ee
  hence
 \begin{eqnarray}
  NT_{j_1\cdots j_n}&=&1-\sum_{i=1}^{n-1}C_{j_1\cdots j_i}\nonumber\\
                    &=&1-\sum_{i=1}^{n-1}NT_{j_1\cdots j_i}\cdot LE_{j_1
                    \cdots j_i}.
 \label{nfs9} \end{eqnarray}                  
  From $NT_{j_1}=1$ and (\ref{nfs9}), we can determine 
  $NT_{j_1\cdots j_n}$. Then from (\ref{nfs6})
  we can determine $C_{j_1\cdots j_n}$ for any generating set 
  $E_{j_1\cdots j_n}$.

  \subsection{  Distribution functions for case 2.}

  We assume that the density function of $\delta_i$ ($i=1,2,\cdots,K$) are $f_i(x)$.
Denote $F_{j_1\cdots j_n}(x)$ the density function of noise term $N_{j_1\cdots j_n}$.
Since the density function of $\xi_{j_1}\delta_{j_1}$ is
$$F_{j_1}(x)=\frac{1}{\xi_{j_1}}f_{j_1}(\frac{x}{\xi_{j_1}}),$$
we have
$$LE_{j_1}=\int_{-\infty}^{-\xi_{j_1}}\frac{1}{\xi_{j_1}}f_{j_1}(\frac{x}{\xi_{j_1}})dx,$$
and $NT_{j_1}=1$, then $C_{j_1}=LE_{j_1}$.
For next stage, since $N_{j_1j_2}=\xi_{j_2}N_{j_1}+\xi_{j_2}\delta_{j_2}$, $N_{j_1j_2}$
will have a density function $F_{j_1j_2}(x)$that is the convolution of the density
functions of $\xi_{j_2}N_{j_1}$ and $\xi_{j_2}\delta_{j_2}$.
If the iteration of $E_{j_1}$ is not terminated we know that
$N_{j_1}$ must greater than $-\xi_{j_1}$, this means the density function of $N_{j_1}$
must equal zero in $(-\infty,-\xi_{j_1})$. Hence, when we take the convolution
 we should use 
 $${\cal F}_{j_1}(x)=\frac{I(x+\xi_{j_1})F_{j_1}(x)}{1-C_{j_1}}$$
 for the density function of $N_{j_1}$, where
$$I(x)=\left\{\begin{array}{ll} 0, &\qquad (x\leq 0)\\
1, &\qquad (x>0) \end{array} \right.$$
We divide by $1-C_{j_1}$ in order to normalize ${\cal F}_{j_1}(x)$. Hence
$$F_{j_1j_2}(x)=\int_{-\infty}^{\infty}\frac{1}{\xi_{j_2}^2}{\cal F}_{j_1}(\frac{t}{\xi_{j_2}})
       f_{j_2}(\frac{x-t}{\xi_{j_2}})dt.$$
$$LE_{j_1j_2}=\int_{-\infty}^{-\xi_{j_1}\xi_{j_2}}F_{j_1j_2}(x)dx,$$
and $NT_{j_1j_2}=1-C_{j_1}$, and hence $C_{j_1j_2}=NT_{j_1j_2}\cdot LE_{j_1j_2}$ is obtained.
  For any generating set $E_{j_1\cdots j_n}$, from (\ref{nfs5}) we will have
  $$LE_{j_1\cdots j_n}=\int_{-\infty}^{-\xi_{j_1}\cdots \xi_{j_n}}F_{j_1\cdots j_n}(x)dx,$$
  where
  $$F_{j_1\cdots j_n}=\int_{-\infty}^{\infty}\frac{1}{\xi_{j_n}^2}{\cal F}_{j_1\cdots j_{n-1}}(\frac{t}{\xi_{j_n}})
  f_{j_n}(\frac{x-t}{\xi_{j_n}})dt,$$
  and $${\cal F}_{j_1\cdots j_{n-1}}(x)=\frac{I(x+\xi_{j_1}\cdots\xi_{j_{n-1}})F_{j_1\cdots j_{n-1}}(x)}{1-C_{j_1\cdots j_{n-1}}}.$$
  Then $NT_{j_1\cdots j_n}=1-\sum_{i=1}^{n-1}C_{j_1\cdots j_i}$, hence we can obtain
  $C_{j_1\cdots j_n}$.

\section
{ Chaotic disturbance on self-similar sets.}

   For any infinite sequence $\{j_i\in I\}_{i=1}^{\infty}$, in this section, 
   in (\ref{nfs2}) $\delta_{j_n}$
   is assigned by a chaotic map. We will first concentrate on the case where
   the chaotic map is a tent map. The tent map is given by
  \be x_{n+1}=\left\{\begin{array}{cc} 2x_n,&\qquad x_n<1/2 \\
   2(1-x_n), &\qquad x_n\geq 1/2.\end{array} \right. \label{nfs10} \ee
 This map is iterated together with the rule given by (\ref{nfs2}). The sequence begins
 with a position $x_0$ chosen arbitrarily. A sequence $\{x_n\}$ is then 
 generated according to (\ref{nfs10})
 and $\delta_{j_n}$ is assigned by
 \be \delta_{j_n}=\left\{\begin{array}{cc} -\epsilon,&\qquad x_n<1/2 \\
   \epsilon, &\qquad x_n\geq 1/2.\end{array} \right. \label{nfs11} \ee
where $\epsilon$ is positive constants less than $1$.
For $\xi$ and $\epsilon$, we assume that
\be \xi+\frac{\epsilon}{1-\xi}<1. \label{nfs12} \ee

  In this section, for any $x_0\in [0,1]$ and any infinite sequence  $\{j_i\in I\}_{i=1}^{\infty}$,
  we will show
 that
 \begin{Theorem}  Under the condition (\ref{nfs12}), the iteration with the rules
   given by (\ref{nfs2}) and (\ref{nfs11}) according to the infinite sequence must
  be terminated at finite order, i.e. there is a generating set at finite
  stage collapse.  From the arbitrarity of the infinite sequence, we can see that
 the self-similar structure is truncated at finite order of iteration.
\end{Theorem}

     First we establish the fact that for $\delta_{j_n}$ given by (\ref{nfs11}), there
  exist an interval $(0,a)$ and $n_0$ such that for $x_0\in (0,a),\ N_{j_1\cdots j_{n_0}}
  <-\xi_{j_1}\cdots\xi_{j_{n_0}}$. If $\delta_{j_i}=-\epsilon,\ (i=1,2,\cdots,n_0)$,
  from (\ref{nfs4}), $N_{j_1\cdots j_{n_0}}
  <-\xi_{j_1}\cdots\xi_{j_{n_0}}$ implies
   $$1+\frac{1}{\xi_{j_1}}+\frac{1}{\xi_{j_1}\xi_{j_2}}+\cdots+\frac{1}{\xi_{j_1}\cdots\xi_{j_{n_0}}}>\frac{1}{\epsilon}.$$
It is sufficient that
$$1+\frac{1}{\xi}+\cdots +(\frac{1}{\xi})^{n_0-1}>\frac{1}{\epsilon},$$
it becomes
$$(\frac{1}{\xi})^{n_0}>\frac{\frac{1}{\xi}-1}{\epsilon}+1,$$
hence
$$n_0>\log(1+\frac{\frac{1}{\xi}-1}{\epsilon})/ \log (\frac{1}{\xi}),$$
then it is sufficient to take $$n_0=[\log(1+\frac{\frac{1}{\xi}-1}{\epsilon})/ \log (\frac{1}{\xi})]+1,$$
where $[A]$ means the integer part of $A$. By knowing $n_0$, for $\delta_{j_i}=
-\epsilon, \ (i=1,2,\cdots,n_0)$ we have
$2^{n_0}a=1/2$ which implies $a=1/2^{n_0+1}$.

  For any $x_0\in (0,1)$, by the ergodicity of the tent map$^{\cite{Sch}}$, after a finite
 number of iterations, falls into the interval $(0,a)$. Suppose it takes $k$ iterations
 to move into $(0,a)$, we have $x_n\in (0,1/2)$ for $k<n<k+n_0$. The noise term
 generated by $x_n$ is
\be N_{j_1\cdots j_n}=-\xi_{j_n}\epsilon-\cdots-\xi_{j_n}\cdots\xi_{j_{k+1}}
\epsilon+
\xi_{j_n}\cdots \xi_{j_{k+1}}N_{j_1\cdots j_k}, \label{nfs13} \ee
where $N_{j_1\cdots j_k}$ can either be positive or negative (when $N_{j_1\cdots j_k}<0$,
 we can assume $N_{j_1\cdots j_k}>-\xi_{j_1}\cdots\xi_{j_k}$, otherwise, we have
 truncated already). We want
to find an integer $l<n_0$ such that
\be N_{j_1\cdots j_{k+l}}<-\xi_{j_1}\cdots\xi_{j_{k+l}}.\label{nfs14} \ee
From (\ref{nfs13}), (\ref{nfs14}) implies
$$(\frac{1}{\xi_{j_{k+l-1}}\cdots\xi_{j_1}}+\cdots+\frac{1}{\xi_{j_k}\cdots\xi_{j_1}})\epsilon-\frac{N_{j_1\cdots j_k}}{\xi_{j_k}\cdots\xi_{j_1}}>1.$$
It is sufficient that
$$\frac{\epsilon}{\xi_{j_1}\cdots\xi_{j_k}}(1+\frac{1}{\xi} +\cdots+(\frac{1}{\xi})^{l-1})-\frac{N_{j_1\cdots j_k}}{\xi_{j_1}\cdots\xi_{j_k}}>1,$$
i.e
$$l>\log(1+\frac{1/ \xi-1}{\epsilon}(\xi_{j_1}\cdots\xi_{j_k}+N_{j_1\cdots j_k}))/ \log(1/ \xi)>0.$$
Then it is sufficient to take
$$l=[\log(1+\frac{1/ \xi-1}{\epsilon}(\xi_{j_1}\cdots\xi_{j_k}+N_{j_1\cdots j_k}))/ \log(1/ \xi)]+1.$$
We must to show $l<n_0$. When $N_{j_1\cdots j_k}\leq 0$,
since $0<\xi_{j_1}\cdots\xi_{j_k}-|N_{j_1\cdots j_k}|<1$, then
$$\log(1+\frac{1/ \xi-1}{\epsilon}(\xi_{j_1}\cdots\xi_{j_k}-|N_{j_1\cdots j_k}|))/ \log(1+\frac{1/ \xi-1}{\epsilon})<1,$$
we obtain $l<n_0$.
When $N_{j_1\cdots j_k}>0$, from (\ref{nfs12}), since 
\begin{eqnarray*}
\xi_{j_1}\cdots\xi_{j_k}+N_{j_1\cdots j_k}& &<\xi_{j_1}\cdots\xi_{j_k}+\frac{1-\xi^{k+1}}{1-\xi}\epsilon\\
& &<\xi^k+\frac{\epsilon}{1-\xi}<\xi+\frac{\epsilon}{1-\xi}<1, 
\end{eqnarray*}
then
$$\log(1+\frac{1/ \xi-1}{\epsilon}(\xi_{j_1}\cdots\xi_{j_k}+N_{j_1\cdots j_k}))/ \log(1+\frac{1/ \xi-1}{\epsilon})<1,$$
we have $l<n_0$.
Thus for the chaotic sequence generated by tent map, our conclusion holds.

  {\bf Remark}: If we define a generating set $E_{j_1\cdots j_n}$ '{\it merges}' when $\xi_{j_1}\cdots\xi_{j_n}\leq N_{j_1\cdots j_n}$,
  we change (\ref{nfs11}) to
   $$\delta_{j_n}=\left\{\begin{array}{cc} \epsilon,&\qquad x_n<1/2 \\
   -\epsilon, &\qquad x_n\geq 1/2,\end{array} \right.  $$
then similarly, for $x_0\in (0,1)$, there exists a generating set at finite
stage merges.

  {\it Generalized case}. From the above discussion, for any chaotic map, if
it is ergodic and there exists an interval $I_0$ such that the same negative
value of $\delta_{j_n}$
are assigned to any $x\in I_0$, then the self-similar structure is truncated at
 a finite order of iteration.

\end{document}